\documentclass[12pt,legno]{article}
\usepackage{amsmath}
\usepackage{amsfonts}
\usepackage{amsthm}

\newcommand{\Rset}{\mathbb{R}}
\newcommand{\Cset}{\mathbb{C}}

\begin{document}

\title{A calculus on L\'evy exponents and selfdecomposability on Banach
spaces\footnote{Research funded by a grant MEN Nr 1P03A04629, 2005-2008.}}

\author{Zbigniew J. Jurek}

\date{February 5, 2008.}

\maketitle

\newtheorem{thm}{THEOREM}
\newtheorem{lem}{LEMMA}
\newtheorem{prop}{PROPOSITION}
\newtheorem{cor}{COROLLARY}

\theoremstyle{remark}
\newtheorem{rem}{REMARK}

\begin{quote}

\noindent {\footnotesize \textbf{ABSTRACT.} In infinite dimensional
Banach spaces there is no  complete characterization of the L\'evy
exponents of infinitely divisible probability measures. Here we
propose \emph{a calculus on L\'evy exponents} that is derived from
some random integrals. As a consequence we  prove that \emph{each}
selfdecomposable measure can by factorized as another
selfdecomposable measure and its background driving measure that is
s-selfdecomposable. This complements a result from the paper of
Iksanov-Jurek-Schreiber in the Annals of Probability \textbf{32},
2004.}



\medskip
\medskip
\underline{Abbreviated title}: \emph{A calculus on L\'evy
exponents}
\end{quote}

\newpage
\noindent \textbf{ 1. Introduction.} Recall that a Borel probability
measures $\mu$, on a real separable Banach space $E$, is called
\emph{infinitely divisible} if for each natural number $n$ there
exists a probability measure $\mu_n$ such that $\mu_n^{\ast n}=\mu$;
the class of all infinitely divisible measures will be denoted by
$ID$. It is well-know that their Fourier transforms (\emph{the
L\'evy-Khintchine formulas}) can be written  as follows
\begin{multline}
\hat{\mu}(y)= e^{\Phi(y)}, \ y \in E', \ \ \mbox{and the exponents
$\Phi$ are of the form} \\  \Phi(y)=i<y,a>- \frac{1}{2}<y,Ry>  +
\int_{E \backslash \{0 \}}[e^{i<y,x>}-1-i<y,x>1_B(x)]M(dx),
\end{multline}
where $E'$ denote the dual Banach space, $<.,.>$ is an appropriate
bilinear form between $E'$ and $E$, $a$ is  a \emph{shift vector},
$R$ is a \emph{ covariance operator} corresponding to the Gaussian
part of $\mu$ and $M$ is a \emph{L\'evy spectral measure}. There is
a one-to-one corresponds between $\mu \in ID$ and the triples
$[a,R,M]$ in its L\'evy-Khintchine formula (1); cf. Araujo-Gin\'e
(1980), Chapter 3, Section 6, p. 136. The function $\Phi(y)$ from
(1) is called the \emph{L\'evy exponent} of $\mu$.

\medskip
\begin{rem} (a) If $E$ is a Hilbert space then L\'evy spectral measures $M$
are completely characterized by the integrability condition
$\int_{E}(1 \wedge ||x||^2)M(dx) < \infty$ and Gaussian covariance
operators $R$ coincide with the positive trace-class operators ; cf.
Parthasarathy (1967), Chapter VI, Theorem 4.10.

(b) When $E$ is \emph{an Euclidean space} then L\'evy exponents are
completely characterized as  continuous negative-definite functions;
cf. Cuppens (1975) and Schoenberg's Theorem on p. 80.
\end{rem}
Finally, a \emph{L\'evy process} $Y(t), t \ge 0$, means a continuous
in probability process with stationary and independent increments
and $Y(0)=0$. Without loss of generality we may and do assume that
it has paths in the Skorochod space $D_E[0,\infty)$ of E-valued
\emph{cadlag functions} (i.e., right continuous with left hand
limits). There is a one-to-one correspondence between the class $ID$
and the class of L\'evy processes.

The cadlag paths of a process $Y$ allows us define \emph{random
integrals} of the form $\int_{(a,b]}h(s)Y(r(ds))$ via the formal
formula of integration by parts. Namely,
\begin{multline}
\int_{(a,b]}h(s)Y(r(ds)):= \\
h(b)Y(r(b))- h(a)Y(r(a)) - \int_{(a,b]} Y(r(s))dh(s), \qquad \qquad
\end{multline}
where $h$ is a real valued function of bounded variation and $r$ is
a monotone and right-continuous function. Furthermore, we have
\begin{equation}
\widehat{\mathcal{L}\Big(\int_{(a,b]}h(s)Y(r(ds))\Big)}(y)= \exp
\int_{(a,b]}\log \widehat{\mathcal{L}(Y(1))}(h(s)y)dr(s),
\end{equation}
where $\mathcal{L}(.)$ denotes the probability distribution and
$\hat{\mu}(.)$ denotes the Fourier transform of a measure $\mu$; cf.
Jurek-Vervaat (1983) or Jurek (1985) or Jurek-Mason (1993), Section
3.6, p. 116.

\medskip
\noindent \textbf{2. A calculus on L\'evy exponents.} Let
$\mathcal{E}$ denotes  the totality of all functions $\Phi :E'\to
\Cset$ appearing as the exponent in the L\'evy-Khintchine formula
(1). Hence we have that
\begin{equation}
\mathcal{E} +\mathcal{E} \subset \mathcal{E} , \ \  \
\lambda\cdot\mathcal{E}\subset \mathcal{E}, \ \ \mbox{for all
postive} \ \ \lambda,
\end{equation}
which means that $\mathcal{E}$ forms a cone in the space of all
complex valued functions defined on $E'$. Furthermore, if $\Phi \in
\mathcal{E}$ then all dilations $\Phi(a \cdot) \in \mathcal{E}$.
These follow from the fact that infinite divisibility is preserved
under convolution and under (convolution) powers to positive real
numbers.

\noindent Here we consider two integral operators acting on
$\mathcal{E}$ or its part. Namely,
\begin{multline}
\mathcal{J}: \mathcal{E} \to \mathcal{E}, \ \
(\mathcal{J}\Phi)(y):=\int_0^1\Phi(sy)ds, \ \ y\in E'; \\
\mathcal{I}: \mathcal{E}_{\log} \to \mathcal{E}, \quad \quad
(\mathcal{I}\Phi)(y):=\int_0^1\Phi(sy)s^{-1}ds, \ \ y\in E'.\qquad
\qquad \qquad \qquad
\end{multline}
Note that $\mathcal{J}$ is well defined on all of$\mathcal{E}$,
since by (3), $\mathcal{J}\Phi$ is the L\'evy exponent of the
well-defined integral $\int_{(0,1]}tdY(t)$, where $Y(1)$ has the
L\'evy exponent $\Phi$; cf. Jurek (1985) or (2004). On the other
hand, $\mathcal{I}$ is only defined on $\mathcal{E}_{\log}$, which
corresponds to infinitely divisible measures with finite logarithmic
moments, since $\mathcal{I}\Phi$ is the L\'evy exponent of the
random integral $\int_{(0,1]}tdY(-\ln t)=\int_{(0,\infty)}
e^{-s}dY(s)$, where $\Phi$ is the L\'evy exponent of $Y(1)$ that has
finite logarithmic moment; cf. Jurek-Vervaat (1983).

\medskip
Here are the main algebraic properties of  the mappings
$\mathcal{J}$ and $\mathcal{I}$.
\begin{lem}
The operators $\mathcal{I}$ and $\mathcal{J}$ acting on appropriate
domains (L\'evy exponents) have the following basic properties:
\begin{description}
\item[(a)] $\mathcal{I,J}$ \ \ \mbox{are additive and positive
homogeneous operators}; \\
\item[(b)] $\mathcal{I,J}$ \ \ \mbox{commute under the composition
and} \ $\mathcal{J}(\mathcal{I}(\Phi))=(\mathcal{I-J})\Phi$.

Other equivalent forms of that last property are: \\
\ \ \ \ \ \ \ \  \ \ \ \ \ \ \ \ \
$\mathcal{J}(I+\mathcal{I})=\mathcal{I}$; \ \
$\mathcal{I}(I-\mathcal{J})=\mathcal{J}$; \ \ \
$(I-\mathcal{J})(I+\mathcal{I})=I$.
\end{description}
\end{lem}
\begin{proof}
Part (a) follows from the fact that $\mathcal{E}$ forms a cone. For
part (b) note that
\begin{multline*}
(\mathcal{J}(\mathcal{I}(\Phi)))(y)=\int_0^1(\mathcal{I}(\Phi))(ty)\,dt=\int_0^1\int_0^1\Phi(sty)s^{-1}dsdt=\\
\int_0^1\int_0^t\Phi(ry)r^{-1}drdt=
\int_0^1\int_r^1\Phi(ry)dt\,r^{-1}dr=\\
\int_0^1\Phi(ry)r^{-1}dr-\int_0^1\Phi(ry)dr=\mathcal{I}\Phi(y)-\mathcal{J}\Phi(y)=
(\mathcal{I}-\mathcal{J})\Phi(y),
\end{multline*}
which proves the equality in (b). Note that from the above (the
first line of  the above argument) we infer also that that operators
$\mathcal{I}$ and $\mathcal{J}$ commute which completes the
argument.
\end{proof}
\begin{lem}
The operators $\mathcal{I}$ and $\mathcal{J}$, defined by (5), have
the following additional properties:
\begin{description}
\item[(a)]$\mathcal{J}:\mathcal{E}_{\log}\to
\mathcal{E}_{\log}$\ \mbox{and} \
$\mathcal{I}:\mathcal{E}_{(\log)^2}\to \mathcal{E}_{\log}$, \\
\item[(b)] $\mbox{If}\ \ (I-\mathcal{J})\Phi \in  \mathcal{E}$ \ \
\mbox{then the corresponding infinitely divisible} \\
\mbox{measure $\tilde{\mu}$ with the L\'evy exponent
$(I-\mathcal{J})\Phi(y)$, $y\in E'$, has finite} \\
\mbox{logarithmic moment}.\\
\item[(c)]$(I-\mathcal{J})\Phi +
\mathcal{I}(I-\mathcal{J})\Phi=(I-\mathcal{J})\Phi +
\mathcal{J}\Phi= \Phi$ \ \ \mbox{for all}\ \ $\Phi\in\mathcal{E}$.
\qquad \qquad \qquad \qquad \qquad
\end{description}
\end{lem}
\begin{proof}
(a) Since the function $E\ni x\to\log(1+||x||)$ is sub-additive, for
an infinitely divisible probability measure $\mu=[a,R, M]$ we have
\begin{multline}
\int_E\log(1+||x||)\mu(dx)<\infty \ \ \mbox{iff}\ \
\int_{\{||x||>1\}}\log(1+||x||)M(dx)<\infty \\
\mbox{iff}\int_{\{||x||>1\}}\log||x||M(dx)<\infty;\qquad \qquad
\qquad \qquad \qquad
\end{multline}
cf. Jurek and Mason (1993), Proposition 1.8.13. Furthermore, if $M$
is the spectral L\'evy measure appearing in  the L\'evy exponent
$\Phi$ then $\mathcal{J}\Phi$ has L\'evy spectral measure
$\mathcal{J}M$ (we keep that potentially conflicting notation),
where
\begin{equation}
(\mathcal{J}M)(A):=\int_{(0,1)}M(t^{-1}A)dt=\int_{(0,1)}\int_E1_A(tx)M(dx)dt,
\end{equation}
for all Borel subsets $A$ of $E\setminus\{0\}$. Hence
\begin{multline*}
\int_{\{||x||>1\}}\log||x||(\mathcal{J}M)(dx)=
\int_{(0,1)}\int_E\,1_{\{||x||>1\}}(tx)\log(t||x||)M(dx)dt \\
 =\int_{(0,1)}\int_{\{||x||>t^{-1}\}}\log(t||x||)M(dx)dt=
\int_{\{||x||>1\}}\int_{||x||^{-1}}^1\log(t||x||)dt\,M(dx)\\
=\int_{\{||x||>1\}}||x||^{-1}\int_{1}^{||x||}\log w\, dw \,M(dx)\\
=
\int_{\{||x||>1\}}||x||^{-1}[||x||\log||x||-||x||+1]M(dx)\\
=\int_{\{||x||>1\}}\log||x||M(dx)-\int_{\{||x||>1\}}[1-||x||^{-1}]M(dx).
\end{multline*}
Since the last integral is always finite as we integrate a bounded
function with respect to a finite measure, we get the first part
of (a). For the second one, let us note that
\begin{multline*}
\int_{\{||x||>1\}}\log||x||(\mathcal{I}M)(dx)=\int_0^{\infty}\int_{\{||x||>1\}}\log||x||M(e^tdx)dt
\\ = 1/2\int_{\{||x||>1\}}\log^2||x||M(dx), \qquad \qquad
\end{multline*}
where $\mathcal{I}M$ is the L\'evy spectral measure corresponding
to the L\'evy exponent $\mathcal{I}\Phi$.

For the part (b), note that the assumption made there implies that
the measure
\begin{equation}
\widetilde{M}(A):=M(A)-\int_{(0,1)}M(t^{-1}A)dt\ge 0, \ \
\mbox{for all Borel sets} \ \ A\subset E\setminus{\{0\}},
\end{equation}
is the L\'evy spectral measure of some $\tilde{\mu}$. [Note that
there is no restriction on the Gaussian part.] In fact, if
$\widetilde{M}$ is a nonnegative measure then it is necessarily a
L\'evy spectral measure because $0\le\widetilde{M}\le M$ and $M$ is
L\'evy spectral measure; comp. Arujo-Gin\'e (1980), Chapter 3,
Theorem 4.7 , p. 119. \noindent To establish the logarithmic moment
of $\tilde{\mu}$ we argue as follows. Observe that for any constant
$k>1$ we have
\begin{multline*}
0\le \int_{(\{1<||x||\le k \}}\log||x||\widetilde{M}(dx)=\\
\int_{\{1<||x||\le k\}}\log||x||M(dx)-\int_{(0,1)}\int_{\{1<||x||\le
k\}}\log||x||M(t^{-1}dx)dt=\\
\int_{\{1<||x||\le
k\}}\log||x||M(dx)-\int_{(0,1)}\int_{\{t^{-1}<||x||\le
kt^{-1}\}}\log(t||x||)\,dM(dx)dt=\\
\int_{\{1<||x||\le k\}}\log||x||M(dx)-\int_{\{1<||x||\le
k\}}\int_{||x||^{-1}}^1\log(t||x||)dt\,M(dx)\\
-\int_{\{k<||x||\}}\int_{||x||^{-1}}^{k||x||^{-1}}\log(t||x||)dt\,M(dx)=\\
\int_{\{1<||x||\le k\}}\log||x||M(dx)-\int_{\{1<||x||\le
k\}}||x||^{-1}\int_1^{||x||}\log(w)\,dw\,M(dx)\\
-\int_{\{k<||x||\}}||x||^{-1}\int_{1}^k\log(w)\,dw\,M(dx)=\\
\int_{\{1<||x||\le k\}}\log||x||M(dx)- \int_{\{1<||x||\le k\}}||x||^{-1}(||x||\log||x||-||x||+1)M(dx) \\
- (k\log k-k+1)\int_{\{||x||>k\}}||x||^{-1}M(dx)= \\
\int_{\{1<||x||\le k\}} (1-||x||^{-1})M(dx)-(k \log k-k+1)\int_{\{||x||>k\}}||x||^{-1}M(dx)\\
\le M(||x||>1)<\infty, \qquad \qquad \qquad
\end{multline*}
and consequently $\int_{(||x||>1)}\log||x||\widetilde{M}(dx<\infty$.
This with property (6), completes the proof of the part (b).

Finally, since $(I-\mathcal{J})\Phi$ is in a domain of definition of
the operator $\mathcal{I}$, so the part (c) is a consequence of
Lemma 1(e) and (d). Thus the proof is complete.
\end{proof}

\medskip
\noindent \textbf{3. Factorizations of selfdecomposable
distributions.} The classes of limit laws $\mathcal{U}$ and $L$ are
obtained by non-linear shrinking transformations and linear
transformations (multiplications by scalars), respectively; cf.
Jurek (1985) and references there. However, there are many
(unexpected) relations between $\mathcal{U}$ and $L$ as was already
proved in Jurek (1985) and more recently in Iksanov-Jurek-Schreiber
(2004). Furthermore, more recently selfdecomposable distributions
are used in modelling real phenomena, in particular in mathematical
finance; for instance cf. Bingham (2006), Carr-Geman-Madan-Yor
(2005) or Eberlein-Keller (1995). This motivates  further studies on
factorizations and other relations between the classes $\mathcal{U}$
and $L$, like those in Theorems 1 and 2, below.

In this section we will apply the operators $\mathcal{I}$ and
$\mathcal{J}$ to L\'evy exponents of selfdecomposable (the class
$L$) and s-selfdecomposable (the class $\mathcal{U}$) probability
measures. For the convenience of the readers recall here that
\begin{multline}
\mu\in L \ \ \mbox{iff}\ \ \forall (t>0) \exists \nu_t \ \
\mu=T_{e^{-t}}\,\mu\ast\nu_t \ \ \\
\mbox{iff} \ \ \mu=\mathcal{L}(\int_{(0,\infty)}e^{-t}dY(t)); \ \
\mathcal{L}(Y(1))\in ID_{\log}, \\
\mu\in \mathcal{U} \ \ \mbox{iff}\ \
 \mu=\mathcal{L}(\int_{(0,1]}\,t\,dY(t)), \ \ \mathcal{L}(Y(1))\in
 ID. \qquad \qquad \qquad \qquad
\end{multline}
Meaures from the class $\mathcal{U}$ are called
\emph{s-selfdecomposable}; cf Jurek (1985), (2004). The
corresponding Fourier transforms of measures from $L$ and
$\mathcal{U}$ easily follow from (2) and (3); cf. Jurek-Vervaat
(1983) or the above references.

\begin{lem}
If $\mu$ is a selfdecomposable probability measure on a Banach space
$E$ with  characteristic function $\hat{\mu}(y)=\exp[\Phi(y)]\, y\in
E',$ then
\[
\widetilde{\Phi}(y):=\Phi(y)-\int_{(0,1)}\Phi(sy)ds=
(I-\mathcal{J})\Phi(y),\ y \in E',
\]
is a L\'evy exponent corresponding to an infinitely divisible
probability measure with finite logarithmic moment.

Equivalently, if $M$ is the L\'evy spectral measure of a
selfdecomposable $\mu$ then the measure $\widetilde{M}$ given by
\[
\widetilde{M}(A):=M(A)- \int_0^1M(t^{-1}A)dt, \ \ A \subset
E\setminus{\{0\}},
\]
is a L\'evy spectral measure on E that additionally integrates the
logarithmic function on the complement of any neighborhood of zero.
\end{lem}
\begin{proof}
If $\mu=[a,R,M]$ is selfdecomposable (or in other words a class L
distribution) then we infer that
\[
M(A)- M(e^tA)\ge 0, \ \ \mbox{for all} \ \ t>0 \ \mbox{and Borel} \
\ A \subset E\setminus{\{0\}},
\]
and that there is no restriction on the remaining two parameters
(the shift vector and the Gaussian covariance operator) in the
L\'evy-Khintchine formula (1). Multiplying both sides by $e^{-t}$
and then integrating over the positive half-line we conclude that
$\widetilde{M}$, given by (8), is a non-negative measure. Since
$\widetilde{M}\le M$ and $M$ is a L\'evy spectral measure, so is
$\widetilde{M}$; comp. Theorem 4.7 in Chapter 3 of Araujo-Gin\'e
(1980). Finally, our Lemma 2(b) gives the finiteness of the
logarithmic moment. Thus the proof is complete.
\end{proof}

\begin{thm}
For each  selfdecomposable probability measure $\mu$, on a Banach
space $E$, there exists a unique s-selfdecomposable probability
measure $\tilde{\mu}$ with finite logarithmic moment such that
\begin{equation}
\mu=\tilde{\mu}\ast\mathcal{I}(\tilde{\mu}) \ \ \mbox{and}  \ \
\mathcal{J}(\mu)=\mathcal{I}(\tilde{\mu})\ .
\end{equation}
In fact, if $\hat{\mu}(y)=\exp[\Phi(y)]$ then
$(\tilde{\mu}){\hat{}}(y)=\exp[\Phi(y)-\int_{(0,1)}\Phi(ty)dt]$,
$y\in E'$.

In other words, if $\Phi$ is the L\'evy exponent of a
selfdecomposable probability measure then $(I-\mathcal{J})\Phi$ is
the L\'evy exponent of an s-selfdecomposable measure with the finite
logarithmic moment and
\begin{equation}
\Phi=(I-\mathcal{J})\Phi+\mathcal{I}(I-\mathcal{J})\Phi=(I-\mathcal{J})\Phi+\mathcal{J}\Phi.
\end{equation}
\end{thm}
\begin{proof}
Let $\hat{\mu}(y)=\exp[\Phi(y)]\in L$. From the factorization in (9)
(the first line) we infer that  $\Phi_t(y):=\Phi(y)-\Phi(e^{-t}y)$
are L\'evy exponents as well. Hence,
\[
\widetilde{\Phi}(y):=\int_{(0,\infty)}\Phi_t(ty)e^{-t}dt=
\Phi(y)-\int_{(0,\infty)}\Phi(e^{-t}y)e^{-t}dt=
((I-\mathcal{J})\Phi)(y)
\]
is a L\'evy exponent as well, because of Lemma 3. Again by Lemma 3
(or Lemma 2 b)), a probability measure $\tilde{\mu}$ defined by the
Fourier transform
$(\tilde{\mu}){\hat{}}(y)=\exp(I-\mathcal{J})\Phi(y)$ has
logarithmic moment. Consequently, $\mathcal{I}(\tilde{\mu})$ is a
well defined probability measure whose L\'evy exponent is equal to
$\mathcal{I}(I-\mathcal{J})\Phi$. Finally, Lemmas 1(b) and 2(c ???)
give the factorization (10).

Since $\mathcal{I}(\tilde{\mu})\in L$ has the property that
$\tilde{\mu}\ast \mathcal{I}(\tilde{\mu})$ is again in $L$,
therefore Theorem 1 from Iksanov-Jurek-Schreiber(2004) gives that
$\tilde{\mu}\in \mathcal{U}$, i.e., it is a s-selfdecomposable
probability distribution.

\noindent To see the second equality in (11) one should observe that
it is equivalent to equality
$\mathcal{J}\Phi=\mathcal{I}(I-\mathcal{J})\Phi$ that indeed holds
true in view of Lemma 1(d).

Suppose there exists another factorization of the form
$\mu=\rho\ast\mathcal{I}(\rho)$ and let $\Xi(y)$ be the L\'evy
exponent of $\rho$. Then we get that $\Phi(y)=\Xi(y)+
(\mathcal{I}\,\Xi)(y)= (I+\mathcal{I})\,\Xi(y)$. Hence, applying to
both sides $\mathcal{I-J}$ we conclude that
\[
(I-\mathcal{J})\Phi=((I-\mathcal{J})(I+\mathcal{I}))\,\Xi=\Xi,
\]
where the last equality is from Lemma 1(b). This proves the
uniqueness of $\tilde{\mu}$ in the representation (10) and thus the
proof of Theorem 1 is completed.
\end{proof}

\begin{rem}
The factorization (10), in Theorem 1, can be also derived from
previous papers as follows:

\emph{for each selfdecomposable (or class $L$) $\mu$ there exists a
unique $\rho\in ID_{\log}$ such that $\mu=\mathcal{I}(\rho)$;
Jurek-Vervaat (1983). Since $\tilde{\mu}:=\mathcal{J}(\rho)$ is
s-selfdecomposable (class) $\mathcal{U}$) with logarithmic moment
(cf. Jurek (1983)) therefore,
$\mathcal{I}(\tilde{\mu})\ast\tilde{\mu}\in L$ in view of
Iksanov-Jurek-Schreiber (2004). Finally, again by Jurek (1985),
$\mathcal{I}(\tilde{\mu})\ast\tilde{\mu}=\mathcal{J}(\mathcal{I}(\rho)\ast\rho)=\mathcal{I}(\rho)=\mu$,
which gives the decomposition.}

However, the present proof is less involved,  more straightforward
and moreover the result and the proof of finiteness of the
logarithmic moment in Lemma 2 (b) is completely new. Last but not
least, the "calculus" on L\'evy exponents, introduced in this note,
is of an interest in itself.
\end{rem}
\begin{rem}
In the case of Euclidean space $\Rset^d$, using Schoenberg's
Theorem, one gets immediately that $\widetilde{\Phi}$ is a L\'evy
exponent; cf. Cuppens (1975), pp. 80-82.
\end{rem}
Following Iksanov, Jurek and Schreiber (2004), p. 1360, we will say
that a selfdecomposable probability measure $\mu $ has \emph{the
factorization property} if $\mu\ast\mathcal{I}^{-1}(\mu)$ is
selfdecomposable as well. In other words, a class $L$ probability
measure convolved with its background driving probability
distribution is again class $L$ distribution. As in
Iksanov-Jurek-Schreiber (2004), Proposition 1, if $L^f$ denotes the
set of all class $L$ distribution with the factorization property
then
\begin{equation}
L^f=\mathcal{I}(\mathcal{J}(ID_{\log}))=
\mathcal{J}(\mathcal{I}(ID_{\log}))=\mathcal{J}(L)
 \ \mbox{and} \ \
L^f\subset L \subset \mathcal{U},
\end{equation}

\begin{cor}
Each selfdecomposable $\mu$ admits a factorization
$\mu=\nu_1\ast\nu_2$, where $\nu_1$ is an s-selfdecomposable measure
(i.e., $\nu_1\in \mathcal{U}$) and $\nu_2$ is a selfdecomposable one
with the factorization property (i.e., $\nu_2\in L^f$). That is,
besides the inclusion $ L^f\subset L \subset \mathcal{U}$ we also
 have that $L\subset L^f \ast \mathcal{U}$.
\end{cor}
\begin{proof}
Because of (10), $\nu_1:=\tilde{\mu}$ is an s-selfdecomposable
measure. Furthermore, $\nu_2:=\mathcal{I}(\tilde{\mu})\in L$ has the
factorization property, i.e., $\nu_2\in L^f$, which completes  the
proof.
\end{proof}
\textbf{EXAMPLES.} 1) Let $ \Sigma_p$ be a symmetric stable
distribution on a Banach space $E$, with the exponent $p$. Then its
L\'evy exponent, $\Phi_p$, is equal to
$\Phi_p(y)=-\int_S|<y,x>|^p\,m(dx)$, where $m$ is a finite Borel
measure on the unit sphere $S$ of $E$; cf. Samorodnitsky and Taqqu
(1994). Hence $(I-\mathcal{J})\Phi_p(y)=p/(p+1)\Phi_p(y)$, which
means that in Corollary 1, both $\nu_1$ and $\nu_2$ are stable with
the exponent $p$ and measures $m_1:=(p/(p+1))m$ and
$m_2:=(1/(p+1))m$, respectively.

  2) Let $\eta$ denotes the Laplace (double exponential)
distribution on real line $\Rset$; cf. Jurek-Yor (2004). Then its
L\'evy exponent $\Phi_{\eta}$ is equal to
$\Phi_{\eta}(t):=-\log(1+t^2),\, t\in\Rset$. Consequently,
$(I-\mathcal{J})\Phi_{\eta}(t)= 2(\arctan t -t)t^{-1}$ is the L\'evy
exponent of the class $\mathcal{U}$ probability measure $\nu_1$ from
Corollary 1, and $(2t-\arctan t -t\,\log(1+t^2))t^{-1}$ is the
L\'evy exponent of the class $L^f$ measure $\nu_2$ from Corollary 1.

\medskip
Before we formulate the next result we need to recall that, by (9),
the class $\mathcal{U}$ is defined here as
$\mathcal{U}=\mathcal{J}(ID)$. Consequently, by  iteration argument
we can define
\begin{equation}
\mathcal{U}^{<1>}:=\mathcal{U}, \ \
\mathcal{U}^{<k+1>}:=\mathcal{J}(\mathcal{U}^{<k>})=\mathcal{J}^{k+1}(ID),
\ k=1,2,... ;
\end{equation}
cf. Jurek (2004) for other characterization of classes
$\mathcal{U}^{<k>}$. Elements from the semigropus
 $\mathcal{U}^{<k>}$ are called \emph{k-times s-selfdecomposable
probability measures}.
\begin{thm}
Let $n$ be any natural number and $\mu$ be a selfdecopmosable
probability measure. Then there exist k-times s-selfdecomposable
probability measures $\tilde{\mu}_k$, for $k=1,2,...,n$, such that
\begin{equation}
\mu=\tilde{\mu}_1\ast\tilde{\mu}_2\ast...\ast\tilde{\mu}_n\ast\mathcal{I}(\tilde{\mu}_n),\
\,\mathcal{J}^k(\mu)=\mathcal{I}(\tilde{\mu}_k), \ \ k=1,2,...,n.
\end{equation}
In fact, if $\Phi$ is the exponent of $\mu$ then $\tilde{\mu}_k$
has the exponent
$\mathcal{I}^{k-1}(I-\mathcal{J})^k\Phi=(I-\mathcal{J})\mathcal{J}^{k-1}\Phi$
and
\begin{multline}
\Phi=(I-\mathcal{J})\Phi
+(I-\mathcal{J})\mathcal{J}\Phi+...+(I-\mathcal{J})\mathcal{J}^{k-1}\Phi+...+
(I-\mathcal{J})\mathcal{J}^{n-1}\Phi+\mathcal{J}^n\Phi \\
=(I-\mathcal{J}^n)\Phi+\mathcal{J}^n\Phi. \quad \quad \quad
\end{multline}
\end{thm}
\begin{proof}
For $n=1$ the factorization (14) and the formula (15) are true by
Theorem 1, with $\tilde{\mu}_1:=\tilde{\mu}$. Suppose our claim (14)
is true for $n$. Since $\rho:=\mathcal{I}(\tilde{\mu}_n)$ is
selfdecomposable, applying to it Theorem 1, we have that
$\rho=\tilde{\rho}\ast\mathcal{I}(\tilde{\rho})$, where
$\tilde{\rho}$ has the L\'evy exponent
$(I-\mathcal{J})\mathcal{J}^n\Phi=\mathcal{J}^n(I-\mathcal{J})\Phi$
 and thus it corresponds to $(n+1)$-times s-selfdecomposable probability because, by Theorem 1,
$(I-\mathcal{J})\Phi$ is already s-selfdecomposable and then we
apply $n$ times the operator $\mathcal{J}$; compare the definition
(13). Thus the factorization (14) holds for $n+1$, which completes
the proof of the first part of the theorem.

\noindent Similarly, applying inductively decomposition (11), from
Theorem 1 and observing from Lemma 1(b) that  we will get the
formula (14). Thus the proof is complete.
\end{proof}

\medskip
\medskip
\textbf{Acknowledgements.} Author would like  to thank the Referee
whose comments improved the language of the paper.

\medskip
\begin{center}
\textbf{REFERENCES}
\end{center}
\medskip
\noindent
[1] A. Araujo and E. Gine (1980). \emph{The central limit theorem for real and
Banach valued random variables.} John Wiley \& Sons, New York.

\noindent [2] N. H. Bingham (2006). L\'evy processes and
self-decomposability in finance, \emph{Probab. Math. Stat.}, vol. 26
Fasc. 2, pp 367-378.

\noindent [3] P. Carr, H. Geman, D. Madan and M. Yor (2005). Pricing
options on realized variance, \emph{Finance and Stochastics}, vol. 9
no 4 , pp. 453-475.

\noindent [4] R. Cuppens (1975). \emph{Decompoqsition of
multivariate probabilities.} Academic Press, New York.

\noindent[5] E. Eberlein and U. Keller (1995). Hyperbolic
distributions in finance, \emph{Bernoulli} vol. 1 , pp. 281-299.

\noindent [6] A. M. Iksanov, Z. J. Jurek, and  B. M. Schreiber
(2004). A new factorization property of the selfdecomposable
probability measures, \emph{Ann. Probab}. vol. 32, No. 2, pp.
1356-1369.

\noindent [7] Z. J. Jurek (1985). Relations between the
s-selfdecomposable and selfdecomposable measures. \emph{Ann.
Probab.} vol.13, No. 2, pp. 592-608.

\noindent [8] Z. J. Jurek (2004). The random integral representation
hypothesis revisited: new classes of s-selfdecomposable laws. In:
Abstract and Applied Analysis; \emph{Proc. International Conf.
ICAAA} , Hanoi, August 2002, p. 495-514. World Scientific, Hongkong.

\noindent [9] Z. J. Jurek and J. D. Mason (1993).
\emph{Operator-limit distributions in probability theory.} John
Wiley \&Sons, New York.

\noindent [10] Z. J. Jurek and W. Vervaat (1983). An integral
representation for selfdecomposable Banach space valued random
variables, \emph{Z. Wahrscheinlichkeitstheorie verw. Gebiete}, 62,
pp. 247-262.

\noindent [11] Z. J. Jurek and M. Yor (2004). Selfdecomposable laws
associated with hyperbolic functions, \emph{Probab. Math. Stat.} 24,
no.1, pp. 180-190.

\noindent[12] K. R. Parthasarathy (1967).\emph{Probabiliy measures
on metric spaces}. Academic Press, New York and London.

\noindent [13] G. Samorodnitsky and M.S. Taqqu (1994).\emph{Stable
non-gaussian random processes}. Chapman \& Hall, New York.

\medskip
\noindent
Institute of Mathematics \\
University of Wroc\l aw \\
Pl.Grunwaldzki 2/4 \\
50-384 Wroc\l aw, Poland \\
e-mail: zjjurek@math.uni.wroc.pl

\end{document}